\documentclass[11pt]{amsart}
\usepackage{amssymb}
\usepackage{verbatim}
\input xypic
\xyoption{all}
\newcommand{\BZ}{{\mathbb{Z}}}

\newcommand{\BC}{{\mathbb{C}}}
\newcommand{\BF}{{\mathbb{F}}}

\newcommand{\BP}{{\mathbb{P}}}

\newcommand{\gb}{\beta}

\newcommand{\gs}{\sigma}
\newcommand{\gS}{\Sigma}
\newcommand{\gO}{\Omega}
\newcommand{\gom}{\omega}

\newcommand{\ga}{\alpha}
\newcommand{\gt}{\tau}

\newcommand{\gT}{\Theta}

\newcommand{\caL}{{\mathcal{L}}}

\newcommand{\ti}[1]{\tilde{#1}}

\newcommand{\ol}[1]{\overline{#1}}

\newcommand{\Pic}{\mathrm{Pic}}
\newcommand{\Div}{\mathrm{Div}}
\newcommand{\jac}{\mathrm{Jac}}

\newcommand{\prym}{\mathrm{Prym}}
\newcommand{\gal}{\mathrm{Gal}}

\newcommand{\dg}{\mathrm{deg}}

\newcommand{\odd}{\mathrm{odd}}
\newcommand{\even}{\mathrm{even}}

\newcommand{\half}{\frac{1}{2}}

\newcommand{\hra}{\hookrightarrow}
\newcommand{\ra}{\longrightarrow}

\swapnumbers
\theoremstyle{plain}
\newtheorem{lma}{Lemma}[section]
\newtheorem{thm}[lma]{Theorem}
\newtheorem{prp}[lma]{Proposition}
\newtheorem{qst}[lma]{Question}
\newtheorem{cor}[lma]{Corollary}

\theoremstyle{definition}
\newtheorem{dfn}[lma]{Definition}
\newtheorem{rmr}[lma]{Remark}
\newtheorem{ntt}[lma]{Notation}

\newtheorem{dsc}[lma]{}

\begin{document}

\title{An explicit formula for the genus 3 AGM}
\thanks{The author was partially supported by Israel-US BSF grant 1998265.}
\author{D. Lehavi}
\address{Department of Mathematics\\
The Ohio State University \\
100 Math Tower \\
231 West 18th Avenue\\
Columbus, OH 43210-1174, USA}
\email{dlehavi@math.ohio-state.edu}
\date{\today}
\begin{abstract}
Given a smooth non-hyperelliptic curve $C$ of genus 3 and a maximal
isotropic subgroup (w.r.t. the Weil pairing) $L\subset\jac(C)[2]$,
there exists a smooth curve $C'$ s.t. $\jac(C')=\jac(C)/L$. This construction
is symmetric. i.e. if we start with $C'$ and the dual flag on it, we get $C$.
A previous less explicit approach was taken by Donagi and Livn\'e (see \cite{DL}).
The advantage of our construction is that it is explicit enough to describe the isomorphism
$H^0(C,\gO_C)\cong H^0(C',\gO_{C'})$.
\end{abstract}
\maketitle

%
\section{Introduction}\label{Sintro}
%
\begin{dsc}
Gauss's {\em arithmetic geometric mean (agM)} can be viewed as an iterative process of dividing a given elliptic
curve $E$, by an element $\ga\in\Pic(E)[2]$, yielding a new elliptic curve
$E'$ (a full treatment is given in \cite{Cox}). This process can be
generalized to the following question:
\end{dsc}
\begin{qst}
For what genera $g$ does there exist a process that, given a curve $C$ of genus
$g$, and a maximal isotropic subgroup (with respect to the Weil pairing)
$L\subset\Pic(C)[2]$, gives a curve $C'$ s.t. $\jac(C')\cong\jac(C)/L$ ?
\end{qst}
\begin{dsc}
This question was answered affirmatively on the case $g=2$ by Richelot
\cite{Ric} and Humbert \cite{Hum} (for a classical treatment, see \cite{BM},
and for a more modern one, see \cite{DL}).
Donagi and Livn\'e proved in \cite{DL}, that no such process can be given for
$g>3$, and presented such a process for $g=3$.
The Donagi-Livn\'e construction has four ``problems'':
\begin{itemize}
  \item It is not symmetric.
  \item There is a 1-parameter choice for $f$, which is unsatisfactory.
  \item It is not easy to give coordinates to the spaces \& functions involved.
  \item It is not obvious how one can Track the canonical class.
\end{itemize}
Our new agM construction avoids these problems. Thus, tracing our construction
one can obtain integration identities (in the spirit of the original agM).
\end{dsc}
\begin{ntt}
In what follows, we fix a ground field $\BF$ such that $\mathrm{char}\BF\neq 2,3$.
Fix a smooth non-hyperelliptic curve $C$ of genus $3$
Lastly fix a maximal isotropic flag (under the Weil pairing) in $\jac(C)[2]$ denoted by
\[
  \caL=\{L_1\subset L_2\subset L_3\}.
\]
Denote by $\ga,\ga'$ the unique non-zero elements in $L_1,L_2^\perp/L_3$ respectively.
\end{ntt}
\begin{dsc}\label{DagM}
We summarize the new agM construction in the diagram:
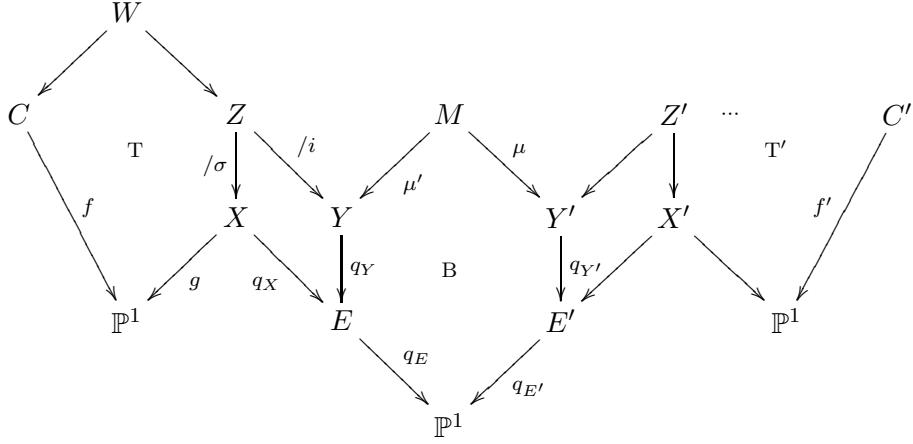
\[
  \xymatrix{
& W \ar [dl] \ar [dr] \\
C \ar [ddr]^{f} \ar@{}[drr]^{\mathrm{T}}& & Z \ar [d]_{/\gs} \ar[dr]^{/i} & & M
\ar[dr]^{\mu} \ar[dl]^{\mu'} \ar@{}[ddd]|{\mathrm{B}}
& & Z' \ar[d] \ar[dl] \ar@{}[r]|{\cdots}
& & C' \ar[ddl]_{f'} \ar@{}[dll]_{\mathrm{T'}}\\
& & X \ar [dl]^{g} \ar[dr]_{q_X} & Y \ar[d]^{q_Y} & & Y' \ar[d]^{q_{Y'}}
& X' \ar[dl] \ar[dr]\\
& \BP^1 & & E \ar[dr]^{q_E} & & E' \ar[dl]^{q_{E'}} & & \BP^1\\
&       & &           & \BP^1 
  }
\]
where:
\begin{itemize}
 \item The zero divisor of the map $f$ is $(f)_0=K_C+\ga$
 \item The zero divisor of the map $f'$ is $(f')_0=K_{C'}+\ga'$
 \item The faces $T,T'$ are trigonal. (see \cite{Do} p. 74)
 \item The construction is symmetric.
 \item The face B is bigonal (see \cite{Do} p. 68-69) with normalization.
 \item The map $q_E:E\ra\BP^1$ identifies two of the branch points of $Y/E$ thus
  describing $L_2$ (see Theorem \ref{Tqis}).
\end{itemize}
The paper is organized as follows:
In section \ref{Sback} we give the necessary technical background, recalling how the double covers
$Y\ra E$ and $F\ra E$ are explicitly determined by the curve $C$ and the element $\ga$, and
vice-versa. In section \ref{Scon} we present the construction, prove the isomorphism:
\[
  \jac(C)/L_3\cong\jac(C'),
\]
and show that the isomorphism
\[
  H^0(C,\gO_C)\cong H^0(C',\gO_{C'})
\]
is explicitly determined by its projectivization
\[
  |K_C|\cong|K_{C'}|.
\]
In section \ref{SKCY} we construct an explicit isomorphism
\[
  |K_C|\cong|K_Y|_\odd
\]
(for the definition of $|K_Y|_\odd$ see Notation \ref{Noddeven}).
In section \ref{SPCR} we show how to construct the double cover $Y'/E'$ from the double cover $Y/E$ and the non zero element in $\ga^\perp/L_2^\perp$, and construct an explicit isomorphism
\[
  |K_Y|_\odd\cong|K_{Y'}|_\odd.
\]
\end{dsc}
\subsection*{Acknowledgments}
The problem was suggested to me by Ron Livn\'e as part of my Ph.D. thesis.
He had positive influence on both the mathematical content, and
the readability of this paper.
%
\section{Technical Background}\label{Sback}
%
We review below the background and necessary results from \cite{Le}.
A proof for all of the assertions in the rest of this section is given there.
\begin{ntt}\label{Nf}
Let $f:C\ra \BP^{1}$ be a function s.t. $(f)_0=K+\ga$.
Denote by $Z,X$ the curves appearing in the reverse trigonal
construction on $f$ (See the \cite{Do} p. 74 and the diagram at \ref{DagM}).
\end{ntt}
\begin{thm}\label{Lthetachar}
The following properties hold:
\begin{enumerate}
  \item There is a natural isomorphism $Z\cong\gT_C\cap(\gT_C+\ga)$.
  \item There are 3 natural commuting involutions on $Z$:
\[
  \begin{aligned}
  i:d&\mapsto d-\ga          \\
  j:d&\mapsto K_C-d          \\
  \gs:=i\circ j:d&\mapsto K_C+\ga-d.
  \end{aligned}
\]
  \item There is a natural isomorphism: $X\cong Z/\gs$.
  \item The fixed points of the involution $j$ are exactly the theta characteristics of $C$ in
    $Z$. Moreover:
    \begin{enumerate}
     \item They are all odd.
     \item They are coupled by $i$ into pairs of the form $z,z+\ga$.
     \item There are 12 such points.
    \end{enumerate}
\end{enumerate}
\end{thm}
\begin{ntt}
Denote
\[
  E:=Z/\langle i,j\rangle,\quad F:=Z/j,\quad Y:=Z/i.
\] 
\end{ntt}
\begin{prp}
The genera of $F,Y,E$ are $1,4,1$ respectively. The curves $E,F$ smooth.
\end{prp}
\begin{thm}\label{Tpc}
The trigonal construction gives a coarse moduli space isomorphism between the moduli of the
following sets of data:
\begin{itemize}
  \item Pairs $C,\ga$ such that the linear system $|K_C+\ga|$ has no base 
  points, up to isomorphisms.
  \item Ramified double covers $Y/E$ of irreducible curves of genera
$4,1$ respectively, such that the curve $E$ is smooth, the curve $Y$ has at most
simple nodes, with a choice of an
element $\gb\in\Pic(E)[2]-\{0\}$, up to isomorphisms.
\end{itemize}
This isomorphism is given by the following maps (which are inverses of one another):
\begin{itemize}
\item The reverse trigonal construction on the map $f:C\ra \BP^1$ (the map $f$ is defined as in
Notation \ref{Nf}) is the double cover $Z\ra X$, which induces the double covers $Y\ra E$ and
$F\ra E$.
\item Given the double covers $Y\ra E,F\ra E$ define $Z:=Y\times_E F$, and define $X$ to be the
third quotient of $Z$, and perform the trigonal construction on one of the rulings of $X$ (they
are conjugate under the involution on $X$).
\end{itemize}
\end{thm}
Theorem \ref{Tpc} can be viewed through the canonical embedding of the curve $C$ (if the curve $C$
is not hyperelliptic then the divisor $|K_C+\ga|$ has no base points):
\begin{thm}\label{CP2prop}
The curve $Z$ can be naturally identified with the set of pairs of points
$\{p_1,p_2\}$ on the curve $C$ such that $p_1+p_2\leq K_C+\ga$.
Under this identification, the following assertions hold:
\begin{enumerate}
  \item The involution $\gs$ takes each such pair to the residual pair of
    points in $K_C+\ga$.
  \item The involution $j$ takes each such pair to the residual 
     pair in $C\cap\ol{p_1 p_2}$. Hence:
  \item The curve $F$ embeds naturally to ${\BP^2}^*$
     as the set of lines through each of these pairs.
  \item The involution $i_F$ acts on the points of $F$ in the
     following way: Let $\{p_1,p_2\}$ be a point in $Z$, and let $\ol{p_1 p_2}$
     be the corresponding line in $F$. There is a unique pair of points
     $\{p_3,p_4\}$ s.t. $K_C+\ga\sim p_1+p_2+p_3+p_4$. We then have 
\[
  i(\ol{p_1 p_2})=\ol{p_3 p_4}
\]
  \item $E\hra\BP^2$ naturally as the Hessian construction (see \cite{Sa} chapter V section V)
    of $i_F$ and the involution induced on the curve $F$ by the involution $i$,
    (as the intersection point of each of the pairs of lines $\{l,i(l)\}$).
  \item There are 6 pairs of bitangents $\{l_{i,1},l_{i,2}\}_{1\leq i\leq 6}$
    such that $\half(l_{i,1}-l_{i,2})=\ga$. These bitangents match (under the canonical embedding)
    the 12 effective theta characteristics in $Z$. (see Lemma \ref{Lthetachar}).
    They are the ramification points of $Z\ra F$.
  \item The 6 intersection points of the pairs $\{l_{i1},l_{i2}\}$ (recall that
   $l_{i2}=i(l_{i1}$) sit on a smooth conic
   $Q\subset\BP^2$ (note that these points are naturally identified with the $q_i$s).
\end{enumerate}
\end{thm}
\begin{ntt}
Choose homogeneous coordinates $x_0,x_1,x_2$ on $\BP^2$.
Denote by $f_E,f_Q\in\BF[x_0;x_1;x_2]$ the homogeneous functions in $x_0,x_1,x_2$ of degrees $3,2$
respectively such that
\[
  \mathrm{Nulls}(f_E)=E,\quad\mathrm{Nulls}(f_Q)=Q.
\]
Denote by $Q_2,Q_3$ respectively the nulls of $f_Q+x_3^2,f_E$ in $\BP^3$ (constructed as a cone
over $\BP^2$). By abuse of notation we denote the nulls of $Q_2,Q_3$ by $Q_2,Q_3$ (if there is
no danger of confusion).
\end{ntt}
\begin{thm}\label{TYcan}
The curve $Y$ is canonically embedded in $\BP^3$ as the complete intersection $Q_2\cap Q_3$.
\end{thm}
\begin{thm}\label{TZiEC}
The trigonal construction and the norm map $\jac Z\ra\jac(Z/i)$
induce an isogeny $\jac(C)\ra\prym((Z/i)/E)$.
The kernel of this isogeny is $\ga^\perp$.
\end{thm}
\begin{ntt}\label{Nqis}
Denote by $q_1,\ldots,q_6$ the points of $\gS/i$
( The $q_i$'s are naturally identified with the 6 ramification points of the double cover
$Y\ra E$).
\end{ntt}
\begin{thm}\label{Tqis}
Considering isotropic (under the Weil pairing) subgroups of $\jac(C)[2]$,
there are bijections between the following objects:
\begin{itemize}
  \item Isotropic groups $L_2\supset\ga$ of order 4, and pairs of $q_i$'s.
  \item Lagrangians $L_3\supset\ga$ and partitions of the $q_i$'s to 3 pairs.
  \item Full flags $L_3\supset L_2\supset\langle\ga\rangle$, and
partitions of the $q_i$s into pairs, with one (of the three) distinguished pair.
\end{itemize}
\end{thm}
\begin{rmr}
The techniques we use here are algebraic. The only requirements we make
are those that Donagi makes in \cite{Do}: $\mathrm{char}\BF\neq 2,3$.
Although the statement of the results in \cite{Le} is for $\BC$, we are using
only the parts that rely on $\mathrm{char}\BF\neq 2,3$.
\end{rmr}
%
\section{The construction}\label{Scon}
%
\begin{ntt}\label{NEf}
Let $q_1,q_2$ be the two $q_i$s (see Notation \ref{Nqis}) that correspond to 
the isotropic subgroup $L_2\subset\jac(C)[2]$ (see Theorem \ref{Tqis}).
Let $\gt:E\ra\BP^1$
be the unique double cover such that $\gt(q_1)=\gt(q_2)$.
Denote by $\ti{H}\ra H$ the bigonal construction on the
tower
\[
  Y\ra E\ra\BP^1.
\]
Denote by $Y',E'$ the normalizations of the curves $\ti{H},H$
respectively.
\end{ntt}
\begin{thm}\label{TAGM}
Assume that the ramification pattern of the tower
\[
  Y\ra E\ra \BP^1
\]
is generic. i.e.
the ramification pattern (see the dictionary at \cite{Do} page 74) is the
following:
\[
  \begin{array}{c|c}
  n(\mathrm{type}) & Y/E \\ \hline
  1 & \subset\subset/= \\
  4 & \subset\subset/\subset \\
  4 & \subset = /= 
  \end{array},
\]
where $n(\mathrm{type})$ is the number of points in $\BP^1$ with a given
non trivial ramification type.
In this case, there exists a curve $C'$ of genus 3, and an element
$\ga'\in\jac(C')[2]$ such that:
\begin{enumerate}
  \item $\jac(C')=\jac(C)/L_3$.
  \item $\langle\ga'\rangle=L_2^\perp/L_3$.
  \item The double cover $Y'\ra E'$ is the double cover related to $C',\ga'$
    in terms of Theorem \ref{Tpc}.
\end{enumerate}
\end{thm}
\begin{proof}
The ramification pattern (see \cite{Do} p. 68-69).
of the bigonal construction on $Y\ra E\ra\BP^1$ is:
\[
  \begin{array}{c|cc}
  n(\mathrm{type}) & Y/E & \ti{H}/H \\ \hline
  1 & \subset\subset/= & \supset\mathrm{\!\!\!\!\!}\subset/\times  \\
  4 & \subset\subset/\subset & \subset = /= \\
  4 & \subset = /= & \subset\subset/\subset
  \end{array},
\]
By \cite{Pa} Proposition 3.1 (page 307) The Abelian variety
$\prym(\ti{H}/H)$ is isomorphic to the dual of the Abelian variety
$\prym(Y/E)$. Since the curve $E$ has only one singular point,
the normalizations of the curves $\ti{H},H$ induce the isomorphism (See \cite{DL} Lemma 1):
\[
  \prym(Y'/E')\cong\prym(\ti{H}/H).
\]
By the Riemann-Hurwitz formula, the genera of the curves $Y',E'$ are
4,1 respectively.
The choice of a partition to two 
pairs on the points of type $\subset =/=$ in the tower $Y\ra E\ra\BP^1$,
is equivalent to a partition to two pairs of the ramification points of
$E'/\BP^1$, which is equivalent to a choice of an unramified cover $F'/E'$.
Define
\[
  Z':=Y'\times_{E'} F'.
\]
By this definition $\gal(Z'/E')=(\BZ/2\BZ)^2$. Define the curve $X'$ to
be third quotient of $Z'$ (the one that is not $F'$ or $Y'$).
By the Riemann-Hurwitz formula, the genera of the curves $Z',X'$ are
7,4 respectively.
Perform the trigonal construction on the double cover $Z'\ra X'$ to get
$f':C'\ra\BP^1$. By Theorem \ref{Tpc} there exists an $\ga''\in\jac(C')[2]$ such
that the zero divisor of the map $f'$ satisfies
\[
  (f')_0=K_{C'}+\ga''.
\]
By Theorem \ref{Tqis} applied to $C',\ga''$ the equality
\[
  \langle\ga''\rangle=L_2^\perp/L_3
\]
holds, thus $\ga''=\ga'$.
\end{proof}
\begin{rmr}
One can describe analogous results for the other (hence, non-generic)
ramification patterns of the tower
\[
  Y \ra E\ra\BP^1.
\]
\end{rmr}
\begin{ntt}
Denote by $i$ the isomorphism 
\[
  i:H^0(C,\gO_C)\ra H^0(C',\gO_{C'}).
\]
Denote the projectivization of this isomorphism by $\BP i$. Denote by $(\BP i)^*$ the induced
isomorphism
\[
  H^*(|K_C|,-)\ra H^*(|K_C'|,\BP i(-)).
\] 
Denote the canonical maps of $C,C'$ by $k_C,k_{C'}$ respectively.
Fix a line $V_\infty\in|K_C|$.
Denote by $z_1,z_2$ affine coordinates on $\BP^2-V_\infty$.
\end{ntt}
\begin{ntt}
Denote by $P_C,P_{C'}$ the maps
\[
  \begin{aligned}
  P_C:H^0(\BP^2,V_\infty)&\ra H^0(C,\gO_C) \\
  l&\mapsto \frac{l dz_1}{\partial k_C/\partial z_2} \\
  P_{C'}:H^0(\BP^2,\BP i(V_\infty))&\ra H^0(C',\gO_{C'}) \\
  l&\mapsto \frac{l \BP i(dz_1)}{\partial k_{C'}/\partial(\BP i(z_2))}
  \end{aligned}
\]
\end{ntt}
\begin{prp}\label{PPC}
The maps $P_C,P_{C'}$ are isomorphisms.
\end{prp}
\begin{proof}
The map
\[
  \begin{aligned}
  H^0(\BP^2,C-3V_\infty)&\ra H^0(\BP^2,V_\infty)\\
  \gom&\mapsto (4V_\infty/k_C)\gom
  \end{aligned}
\]
is an isomorphism since the degree map $\dg:\Pic(\BP^2)\ra\BZ$ is an isomorphism.
Up to this map, the map $P_C$ is the Poincar\'e residue map
\[
  H^0(\BP^2,C-K_{\BP^2})\ra H^0(C,K_C),
\]
which is an isomorphism in this case (see \cite{GH} p. 221).
\end{proof}
\begin{cor}
The equality
\[
  i=P_{C'}\circ(\BP i)^*\circ P_C^{-1}.
\]
holds
\end{cor}
\begin{proof}
This follows from Proposition \ref{PPC}.
\end{proof}
%
\section{Analysis of $|K_C|$ and $|K_{Y}|$}\label{SKCY}
%
\begin{lma}\label{Leo}
Let $\gom\in H^0(Y,\gO_Y)$ be a differential such that the zero divisor
 $\gom_0$ is symmetric under
the involution  $i$. Then the differential $\gom$ is either symmetric or
antisymmetric under the involution $i$.
\end{lma}
\begin{proof}
The zero divisor of the differential $\gom+i\gom$ satisfies
\[
 (\gom+i\gom)_0\supset(\gom)_0.
\]
Hence, there exists some $k\in\BC$ such that
\[
 k(i\gom)=k\gom=\gom+i\gom.
\]
So either $k=0$ or $k=2$.
\end{proof}
\begin{ntt}\label{Noddeven}
Denote the the odd (respectively even) part of the cohomology group
$H^0(Y,\gO_Y)$ under the double cover $Y\ra E$ by
$H^0(\gO_Y,Y)_{\odd}$ (respectively $H^0(\gO_Y,Y)_{\even}$).
Denote by $|(K_Y)_{\odd}$ the projectivization of the vector space $H^0(\gO_Y,Y)_{\odd}$.
Given a cover $A\ra B$ denote the ramification divisor of the cover by
$R_{A/B}$.
\end{ntt}
\begin{rmr}
Recall (see Theorem \ref{TYcan}) that under the canonical map, the curve $Y$
is mapped to the complete intersection $Q_2\cap Q_3$, where $Q_3$
is a cone over $E\subset H$ and $H\subset\BP^3$ is a plane.
\end{rmr}
\begin{lma}\label{LYdeco}
The following properties hold:
\begin{enumerate}
  \item $H^0(Y,\gO_Y)_{\even}$ is one dimensional, and its divisor is $R$.
  \item $\BP(H^0(Y,\gO_Y)_{\even})\subset|K_Y|$ is the point dual to $H$.
  \item $H^0(Y,\gO_Y)_{\odd}$ is $3$ dimensional.
  \item $\BP(H^0(Y,\gO_Y)_{\odd})\subset|K_Y|$ is the dual space
    of the hyperplanes in $\BP^3$ through the vertex of $Q_3$.\label{Iiso}
\end{enumerate}
\end{lma}
\begin{proof}
The first two claims follow from Theorem \ref{TYcan}.
The third claim comes from the decomposition
\[
  H^0(Y,\gO_Y)=H^0(Y,\gO_Y)_{\odd}\oplus H^0(Y,\gO_Y)_{\even}.
\]
We are left with the last claim.
Let $H'\subset\BP^3$ be a plane through the vertex of the cubic surface $Q_3$.
Denote $D:=Y\cap H'$. Since $\dg(Y)=6$, the number of points in the set
$D$ (counting multiplicities) is 6.
Let $\gom\in H^0(Y,\gO_Y)$ be a differential
on the curve $Y$ such that the zero divisor $(\gom)_0$ satisfies
\[
  (\gom)_0=D.
\]
Since the set $D$ is invariant under the involution $\gs$,
Lemma \ref{Leo} implies that the differential $\gom$ is either symmetric or
anti-symmetric.
As the even differentials are analyzed above, the differential $\gom$ is odd.
Since the space of such planes $H'\subset\BP^3$ is 2 dimensional, this space
is all of $\BP(H^0(\gO_Y,Y)_{\odd})^*$.
\end{proof}
\begin{dfn}\label{Dphi}
Denote by $\phi_Y$ the map
\[
  \begin{aligned}
  \phi_Y:H^*&\ra\BP(H^0(Y,\gO_Y)_{\odd}) \\
  L^*&\mapsto {q_Y}^* L\cap E.
  \end{aligned}
\]
\end{dfn}
\begin{cor}\label{Cphiiso}
The map $\phi_Y$ is a natural isomorphism.
\end{cor}
\begin{proof}
This follows from composing the isomorphism of claim \ref{Iiso} in Lemma
\ref{LYdeco} with the natural isomorphism
\[
  \begin{aligned}
  H^*&\ra\text{ planes through the vertex of }Q_3 \\
  L^*&\mapsto\text{ the unique plane through the vertex and }L.
  \end{aligned}
\]
\end{proof}
\begin{prp}\label{Pisotri}
The norms in the trigonal construction on $f$ induce an isomorphism
\[
  H^0(C,\gO_C)\overset{\cong}{\ra}H^0(Z,\gO_Z)/H^0(X,\gO_X).
\]
\end{prp}
\begin{proof}
By \cite{Do} Theorem 2.11 (p. 76) the norms in the trigonal construction induce an
isomorphism:
\[
  \jac(C)\cong\prym(Z/X).
\]
Whence, they induce an isomorphism on the universal covers.
\end{proof}
\begin{ntt}\label{NmCa}
Denote by
 \[
  t:H^0(C,\gO_C)\ra H^0(Y,\gO_Y)_{\odd}
\]
The isomorphism given by the composition of the norms at diagram in \ref{DagM}.
(it is an isomorphism since the map $\jac(C)\ra\prym(Y/E)$ is an isogeny, see Theorem
\ref{TZiEC}).
Denote by $\psi$ the map
\[
  \begin{aligned}
  \psi:H^*&\ra\BP(H^0(C,\gO_C))\\
  L^*&\mapsto L\cap C.
  \end{aligned}
\]
Denote by $m_{C,\ga}$ the automorphism of the plane $H$ dual to the
automorphism
\[
  \phi_Y^{-1}\circ\BP(t)\circ\psi
\]
of the dual plane $H^*$.
\end{ntt}
\begin{prp}
The following properties hold:
\begin{itemize}
  \item The inequality $R_{W/C}>R_{W/Z}$ holds (in $\Div(W)$).
  \item The divisor $R_{W/Z}$ is invariant under the involution $\gt$.
\end{itemize}
\end{prp}
\begin{proof}
The first claim follows from the trigonal construction dictionary (see \cite{Do} p. 74).
The second claim is true since the curve $Z$ is defined
as the quotient of $W$ by the involution $\gt$.
\end{proof}
\begin{thm}
The automorphism $m_{C,\ga}$ is the identity.
\end{thm}
\begin{proof}
We will prove that the automorphism $m_{C,\ga}$ fixes each of the 6 points of
the set $Q\cap E$. Since the set
$Q\cap E$ is not contained in a line, this will prove our claim.
Let $L$ be one of the 12 bitangents $l_{ij}$.
Let $\gom\in H^0(C,\gO_C)$ be a representative of
$\psi(L)\in\BP(H^0(C,\gO_C))$.
Denote the two tangency points of the line $L$ and the curve $C$ by $p,q$.
Denote by $r,t$ the two residual points in the divisor class $|K_C+\ga|$.
We will calculate the zero divisor of the differential:
\[
  \xi:=s_*(\pi^*\gom+s\pi^*\gom)\in H^0(\gO_Z,Z).
\]
By \cite{Ha} proposition IV.2.1, 
\[
  \begin{aligned}
  (\pi^*\gom)_0&=\pi^*((\gom)_0)+R_{W/C}, \\
\text{and }(\pi^*\gom+s\pi^*\gom)_0&=s^*((\xi)_0)+R_{W/Z}.
  \end{aligned}
\]
Whence,
\[
  \begin{aligned}
  \pi^*((\gom)_0)&=\{x,\ti{x}|x\in L\cap C\text{ and }x+\ti{x}\geq K_C+\ga\}\\
                 &=2((p,q)+(p,r)+(p,t)+(q,p)+(q,r)+(q,t)) \qquad&\Rightarrow\\
  (\pi^*\gom)_0&=\pi^*((\gom)_0)+R_{W/C}\geq R_{W/Z}+2((p,q)+(q,p))
&\Rightarrow \\
  (\pi^*\gom+s\pi^*\gom)_0&\geq R_{W/Z}+2((p,q)+(q,p)) &\Rightarrow \\
  (\xi)_0&\geq 2\{p,q\}
\end{aligned}
\]
By Proposition \ref{Pisotri}, the set $(\xi)_0$ is invariant under the
involution $\gs$. Since
\[
  K+\ga-(p+q)=r+t,
\]
we have:
\[
  \begin{aligned}
  (\xi)_0&\supset 2\{p,q\}+2\{r,t\}. \qquad \Rightarrow \\
  ({q_Y}_*(\xi))_0&\supset 2\{\{p,q\},\{r,t\}\}.
  \end{aligned}
\]
The point $\{\{p,q\},\{r,t\}\}$ is one of the ramification points of the double
cover
$Y\ra E$. Under the canonical embedding of the curve $Y$ (see Theorem
\ref{TYcan}) this point is mapped to the point $\ol{pq}\cap\ol{rt}\in H$.
i.e. the line
$L=\ol{pq}$ is mapped under $m_{C,\ga}$ to a line through the point
$\ol{pq}\cap\ol{rt}$.
By symmetry considerations, $\ol{pq}\cap\ol{rt}$ is mapped to itself.
Since this applies to any one of the 6 points in the set $Q\cap E$, we have
proved the claim.
\end{proof}
%
\section{Analysis of $|K_Y|$ and $|K_{Y'}|$}\label{SPCR}
%
\begin{ntt}\label{NKoddeven}
Denote $M:=Y\times_{\BP^1}E'$. Denote by $\mu,\mu'$ the
quotients
\[
  \mu:M\ra Y,\qquad \mu':M\ra Y'.
\]
Denote by $\gt'$ the
involution on $E'$ s.t. $E/\gt'=\BP^1$. Denote by $q_E,q_{E'}$ the quotients by
the involutions $\gt,\gt'$ respectively.
Denote by $q_Y,q_{Y'}$ the quotient by the double covers $Y\ra E,Y'\ra E'$.
Recall Notation \ref{Noddeven}.
We use the analog notations $H^0(\gO_{Y'},Y')_{\odd},H^0(\gO_{Y'},Y')_{\even},|(K_{Y'})_{\odd}|$ for the curve $Y'$.
\end{ntt}
\begin{lma}\label{Lbigoiso}
The norms in the bigonal face of the diagram in \ref{DagM} (face ``B'') induce an
isomorphism
\begin{equation}\label{EYY}
  \mu'_*\circ\mu^*:H^0(Y,\gO_Y)_{\odd}\overset{\cong}{\ra}
                   H^0(Y',\gO_Y')_{\odd}.
\end{equation}
\end{lma}
\begin{proof}
By \cite{Pa}  Proposition 3.1 (p. 307) the norms in the bigonal construction induce an
isomorphism:
\[
  \prym(Y/E)\cong\prym(Y'/E')^\wedge.
\]
Whence, they induce an isomorphism on the universal covers.
\end{proof}
\begin{ntt}
Denote by $P$ the isomorphism 
\[
  P:|(K_Y)_{\odd}|\ra|(K_{Y'})_{\odd}|
\]
induced from the isomorphism in Equation \ref{EYY}.
\end{ntt}
\begin{dsc}
We will describe the isomorphism $P$ with two goals in mind:
\begin{itemize}
\item Combined with the isomorphisms 
\[
  |K_C|\cong|(K_Y)_{\odd}|,\text{ and }|K_C'|\cong|(K_{Y'})_{\odd}|
\]
(described in Notation \ref{NmCa})
we will get an explicit description of the isomorphism
\[
  \BP H^0(C,\gO_C)=|K_C|\cong|K_{C'}|=\BP H^0(C',\gO_{C'}).
\]
\item Using this isomorphism can calculate the plane
configuration of $Y'$, from the plane configuration of $Y$.
\end{itemize}
\end{dsc}
\begin{ntt}\label{Nphi}
We will use the notation
$\phi_Y$, presented in Definition \ref{Dphi}. Define $\phi_{Y'}$ analogy.
\end{ntt}
\begin{dfn}\label{DLp}
Considering the curves $E,E'$ embedded in $\BP^2$, the maps $\gt,\gt'$ are
given as projections from points on
\[
  t\in E\subset\BP^2,\qquad t'\in E'\subset\BP^2
\]
respectively. For any point  $p\in\BP^2$, denote by $L_p$ (respectively
$L_P'$) the line between the two points of $q_E^{-1}(p)$ (respectively 
$q_{E'}^{-1}(p)$).
\end{dfn}
\begin{prp}\label{PLp}
The following equalities hold:
\[
  L_p\cap E=q_E^{-1}(p)+t,\qquad L_p'\cap E'=q_{E'}^{-1}(p)+t.
\]
\end{prp}
\begin{proof}
This follows from Definition \ref{DLp}
\end{proof}
\begin{thm}\label{TLp}
Let $p$ be a point in $\BP^1$, then
\[
  P(\phi_Y(L_p))=\phi_{Y'}(L_p').
\]
\end{thm}
\begin{proof}
Let $\gom\in H^0(Y,\gO_Y)_{\odd}$ be a representative of
\[
  \phi_Y(L_p)\in\BP(H^0(Y,\gO_Y)_{\odd}).
\]
By \cite{Ha} proposition IV.2.1:
\[
  (\mu^*\gom)_0\supset\mu^*((\gom)_0)=\mu^*(q_Y^*(L_p\cap E))\supset
    (q_E\circ q_Y\circ\mu)^{-1}(p).
\]
All The points in the effective divisor $(q_E\circ q_Y\circ\mu)^{-1}(p)$ are
 moving with $p$. Whence, generically,
\[
  R_{M/Y'}\cap(q_E\circ q_Y\circ\mu)^{-1}(p)=\emptyset.
\]
By \cite{Ha} proposition IV.2.1:
\[
  \begin{aligned}
  (\mu^*\gom)_0=&\mu'_*(\mu^*\gom))_0+R_{M/Y'}\qquad\Rightarrow\\
  2(\mu'_*(\mu^*\gom))_0=&\mu'_*((\mu^*\gom)_0)-\mu'_*(R_{M/Y'})
  \supset\mu'_*((q_E\circ q_Y\circ\mu)^{-1}(p)) \\
  =&\mu'_*((q_{E'}\circ q_{Y'}\circ\mu')^{-1}(p))
  =2(q_{E'}\circ q_{Y'})^{-1}(p).
  \end{aligned}
\]
By Lemma \ref{Lbigoiso} the differential $\mu'_*((\mu^*\gom)$ is odd, and by proposition 
\ref{PLp}, the differential
\[
  \mu'_*(\mu^*\gom)\in H^0(Y',\gO_{Y'})_{\odd}
\]
is a representative of $\phi_{Y'}(L_p')$.
\end{proof}
\begin{cor}\label{CLp}
The equality
\[
  P(t)=t'
\]
hold.
\end{cor}
\begin{proof}
This follows from the definition of $t$ and $t'$ (Definition \ref{DLp}).
\end{proof}
\begin{lma}\label{LBMY}
There are no multiple points in the ramification divisor $R_{M/\BP^1}$.
\end{lma}
\begin{proof}
By the bigonal construction dictionary (see \cite{Do} p. 68-69) for any point
$p\in\BP^1$, if the cover $Y'\ra\BP^1$ is ramified, then the double cover
$E\ra\BP^1$ is etale.
The result follows since $M=Y'\times_{\BP^1}E$, and there are no multiple
points in $R_{Y'/\BP^1}$.
\end{proof}
\begin{thm}\label{T4pts}
The following equalities (of sets of points in $\BP^2$) hold:
\[
  P(Q\cap E\smallsetminus\{q_1,q_2\})=B_{E'/\BP^1}, \qquad
  P^{-1}(Q'\cap E'\smallsetminus\{q'_1,q'_2\})=B_{E/\BP^1}.
\]
\end{thm}
\begin{proof}
Let $p\in\BP^1$ be a branch point of the double cover $q_E:E\ra\BP^1$.
Denote by $q$ the ramification point of the map $q_E$ above $p$.
Let $\gom$
be an odd differential on the curve $Y$ such that
$(\gom)_0\supset q_Y^{-1}(q)$.
i.e. the line $\phi_Y^{-1}(\BP(\gom))\subset\BP^2$ (see Definition \ref{Dphi})
passes through the the point $q\in\BP^2$.
By the bigonal construction dictionary (see \cite{Do} p. 68-69) there is a unique
point $q'\in q_{E'}^{-1}(p)$ such that the double $q_{Y'}:Y'\ra E'$ is branched
at $q'$.
By \cite{Ha} proposition IV.2.1:
\[
  \begin{aligned}
  (\mu^*\gom)_0\supset&\mu^*((\gom)_0)\supset\
     \half(q_E\circ q_Y\circ\mu)^{-1}(p)=
     \half(q_{E'}\circ q_{Y'}\circ\mu')^{-1}(p) \\
   \supset&\half(q_{Y'}\circ\mu')^{-1}(q').
  \end{aligned}
\]
However, by Lemma \ref{LBMY}
\[
  (q_{Y'}\circ\mu')^{-1}(q')\cap R_{M/Y'}=\emptyset.
\]
By Lemma \ref{Lbigoiso},
\[
  2\mu^*\gom={\mu'}^*\mu'_*\mu^*\gom.
\]
Whence, by \cite{Ha} proposition IV.2.1:
\[
  \begin{aligned}
  (\mu^*\gom)_0=&{\mu'}^*(\mu'_*(\mu^*\gom)_0)+R_{M/Y'}\qquad\Rightarrow\\
  2(\mu'_*(\mu^*\gom))_0=&\mu'_*((\mu^*\gom)_0)-\mu'_*(R_{M/Y'})
  \supset\half\mu'_*((q_{Y'}\circ\mu)^{-1}(q')) \\
  =&q_{Y'}^{-1}(q').
  \end{aligned}
\]
By Lemma \ref{Lbigoiso} the differential $\mu'_*(\mu^*\gom)\in H^0(Y',\gO_{Y'})$ is odd.
Whence, the point $q'\in\BP^2$ sits on the line
\[
 \phi_{Y'}^{-1}(\mu'_*(\mu^*\gom))
\]
(see Notation \ref{Nphi} and Definition \ref{Dphi} for the definition of
$\phi_{Y'}$).
Thus:
\[
  P(q)=q'.
\]
The second assertion is symmetric.
\end{proof}
\begin{dsc}
Corollary  \ref{TLp} and Theorem \ref{T4pts} enable us to calculate directly
the plane
configuration of the double cover $Y'\ra E'$ from the plane configuration of
the double cover $Y'\ra E'$.
The curve $E'\subset\BP^2$ is the unique cubic that passes through the
following 9 points
\begin{itemize}
\item The point $t$ (see Definition \ref{DLp}).
\item The four ramification points of the double cover $q_E:E\ra\BP^1$.
\item The four points of $Q\cap E\smallsetminus\{q_1,q_2\}$.
\end{itemize}
such that the lines through the point $t$ are tangent to $E'$ at the four
points of $Q\cap E\smallsetminus\{q_1,q_2\}$.
The conic $Q'\subset\BP^2$ is the unique conic that passes through the
following 6 (note that there are 6, and not only 5 known points):
\begin{itemize}
\item The four ramification points of the double cover $q_E:E\ra\BP^1$.
\item The two points of $\ol{q_1q_2}\cap E'$.
\end{itemize}
\end{dsc}

\end{document}